\documentclass[11pt, psamsfonts]{amsart}
\usepackage{amsmath}
\usepackage{amsthm}
\usepackage{amssymb}
\usepackage{amscd}
\usepackage{amsfonts}
\usepackage{amsbsy}
\usepackage{epsfig}

\theoremstyle{plain}
\newtheorem{thm}{Theorem}

\newtheorem{lem}[thm]{Lemma}

\theoremstyle{definition}

\newtheorem{rmk}[thm]{Remark}

\newcommand{\st}{such that }

\newcommand{\DDD}{\mathbb{D}}
\newcommand{\NNN}{\mathbb{N}}

\newcommand{\QQQ}{\mathbb{Q}}

\newcommand{\CCC}{\mathbb{C}}

\newcommand{\h}{\mathcal{H}}

\begin{document}

\begin{abstract}
We characterise the set of fixed points of a class of holomorphic
maps on complex manifolds with a prescribed homology.  Our main
tool is the Lefschetz number and the action of maps on the first
homology group.
\end{abstract}

\title[Periods for holomorphic maps via Lefschetz numbers]
{Periods for holomorphic maps\\ via Lefschetz numbers}
\author[Jaume Llibre and Michael Todd]{Jaume Llibre and Michael
Todd}
\address{Departament de Matem\`{a}tiques,
Universitat Aut\`{o}noma de Barcelona,\newline 08193 Bellaterra,
Barcelona, Spain} \email{jllibre@mat.uab.es}

\address{Mathematics Department,
University of Surrey, Guildford, \newline Surrey, GU2 7XH, UK}
\email{m.todd@surrey.ac.uk}

\thanks{The first author was partially supported by a
MCYT grant BFM2002--04236--C02--02 and by a CIRIT grant number
2001SGR 00173, the second one by a Marie Curie Fellowship number
HPMT-CT-2001-00247. \newline
This paper has been accepted by
Abstract and Applied Analysis.}



\maketitle

\section{Introduction and main result}
In this note we are concerned with fixed point theory for
holomorphic self maps on complex manifolds. After the well--known
Schwarz Lemma on the unit disk, which assumes a fixed point, the
Pick Theorem was proved in \cite{pick}.  This can be extended to a
Pick--type theorem on hyperbolic Riemann surfaces as is shown in
\cite{milnonedim} and \cite{hou}. For a more general type of
space: open, connected and bounded subsets of a Banach space, the
Earle-Hamilton Theorem was proved in \cite{earlham}.

We prove a similar theorem for periodic points of holomorphic maps
of complex manifolds with a prescribed homology. Our proof comes
from a topological viewpoint and bears no relation with other
methods of proof for this problem. We use a result of the first
author and Fagella in \cite{holo} which relates the Lefschetz
number to the set of fixed points for holomorphic maps on compact
complex manifolds.  In fact, we show that for many situations,
holomorphic maps can have at most one fixed point.

We recall the classical results first. First we state a Pick--type
theorem.  This is very close to the well--known Pick Theorem on
the unit disk $\DDD$, which is in turn very similar to the
classical Schwarz Lemma on $\DDD$ (see \cite{milnonedim} for both
results). Let ${\rm Per}(f)$ denote the set of periods of the
periodic points of $f$. Also define ${\rm Fix}(f)$ to be the set
of fixed points of $f$ and define ${\rm Per}_m (f)$ to be the set
of periodic points of period $m$.

The following Pick-type theorem is Theorem 5.2 of
\cite{milnonedim}.

\begin{thm}
Suppose that $S$ is a hyperbolic Riemann surface.  Then for a
holomorphic map $f: S \to S$,  \begin{itemize} \item[(a)] either
$f^p$ is the identity for some $p$; \item[(b)] or $f$ has one
periodic point: a fixed point; \item[(c)] or $f$ no periodic
points.\end{itemize}

\label{thm:riem}

\end{thm}

The following is the Earle--Hamilton Theorem, see \cite{earlham}
and \cite{harris} for more details.  (In fact, the usual setting
for this theorem is complex Banach spaces.)   We say that a set
$E$ is mapped {\em strictly inside}  a subset $F$ of the metric
space $X$ by the map $f$ if there is some $\epsilon>0$ \st
$B_\epsilon(f(E)) \subset F$ where $B_\epsilon(A)$ is the
$\epsilon$--neighbourhood of the set $A \subset X$.

\begin{thm}
Let $M \subset \CCC^m$ be a bounded connected open set and $f:M
\to M$ be holomorphic. If $f(M)$ lies strictly inside $M$ then $f$
has a unique fixed point.

\label{thm:earlham}

\end{thm}

Note that if the map $f$ satisfies the theorem then so does $f^m$.
Therefore, ${\rm Per}(f)=\{1\}$.

Given a complex manifold $M \subset \CCC^m$, let $f:M \to M$ be a
holomorphic map which extends continuously to $\overline{M}$ the
closure of $M$, and furthermore, so that $f(\overline{M})$ is
inside $\stackrel{\circ}{M}$, the interior of $M$. Here we take
the closure with respect to the usual topology on $\CCC^m$ and the
interior with respect to the the topology of $M$. We denote the
set of such maps by $\h (M)$.  Observe that if $M$ has no boundary
then whenever $f: M \to M$ is holomorphic, we have $f \in \h (M)$.
Also note that whenever $M$ is open and bounded then $f \in \h
(M)$ maps $M$ strictly inside itself. The following is our main
result.

\begin{thm}
Suppose that $M$ is a complex manifold  and $f \in \h (M)$.
Suppose further that $H_1(M, \QQQ) = \stackrel{n}{\overbrace{\QQQ
\oplus \cdots \oplus\QQQ}}$ for $n \geq 0$ (when $n=0$ we suppose
$H_1(M, \QQQ) =\{0\}$) and $H_k(M, \QQQ) = 0$ for $k>1$. Then
\begin{itemize}
\item[(a)] either there is some $m \geq 1$ \st $f^m$ has infinitely many
fixed
points;
\item[(b)] or $f$ has only one periodic point: a fixed point;
\item[(c)] or $f$ has no periodic points.

\end{itemize}
Furthermore, if $H_k(M, \QQQ) = 0$ for $k>0$ then we are in case
(b). \label{thm:topick}
\end{thm}

For example, in case (a) we could consider a rational rotation of
the torus; in case (b) we could consider the map $z \mapsto
\frac{z}{2}$ on the unit disk; in case (c) we could consider an
irrational rotation of the torus.

Note that our theorem can apply in any dimension and does not
require manifolds to be open or bounded.  Therefore, in some cases
it is an extension of Theorems~\ref{thm:riem} and
\ref{thm:earlham}.

\section{Proof of the main result}
\label{sec:complex}

To prove Theorem~\ref{thm:topick} we use Lefshetz number and
Theorem~\ref{thm:holo} below. We recall the concept of Lefschetz
number. Let $M$ be a compact manifold of dimension $M$. Then a
continuous map $f:M \to M$ induces an endomorphism $f_{\ast k}:
H_k(M, \QQQ) \to H_k(M, \QQQ)$ for $k=0, 1, \ldots, n$ on the
rational homology of $M$. The {\em Lefschetz number of $f$} is
defined by $$L(f) = \sum_{k=0}^n (-1)^k {\rm trace}(f_{\ast k}).$$

Since $f_{\ast k}$ are integral matrices, it follows that $L(f)$
is an integer.

By the well known Lefschetz Fixed Point Theorem, if $L(f) \neq 0$
then $f$ has a fixed point (see, for instance, \cite{brown}). We
can consider $L(f^m)$ too: also $L(f^m) \neq 0$ implies that $f^m$
has a fixed point.

Our main tool here is the following theorem of \cite{holo}.

\begin{thm}
Let $M$ be a complex manifold and $f \in \h (M)$ be a nonconstant
map.  Suppose that for $m \geq 1$, all fixed points of $f^m$ are
isolated. Then $L(f^m) \geq \# {\rm Fix}(f^m)$. \label{thm:holo}

\end{thm}

In fact, in \cite{holo} the assumption on the periodic points was
that the set of all periodic points must be isolated.  It is
straightforward to weaken this assumption as above.  Also, it was
assumed that $f$ was holomorphic on a compact manifold $M$ and
that $f(M) \subset \stackrel{\circ}{M}$.  However, we don't need
holomorphicity on the boundary of our manifold since we can't have
any fixed points there.

\begin{rmk}
In fact, it is easy to see from the proof of this theorem, that we
come to the same conclusions for maps which are not holomorphic,
but for which det$(I-Df^m(x))>0$ for any $x \in {\rm Fix}(f^m)$.
\label{rmk:posindex}
\end{rmk}

We will need the following lemma, see the appendix of
\cite{milngen}.

\begin{lem}
Let $\mu_1, \ldots, \mu_q$ be complex numbers with $|\mu_j| = 1$.
For any $\epsilon>0$ there exist infinitely many values of $m$ \st
$|\mu_j^m -1|< \epsilon$ for all $1\leq j \leq q$. Hence,
$Re(\mu_j^m)>1-\epsilon$ for every $\mu_j$.

\label{lem:modone}

\end{lem}

{\it Proof of Theorem~\ref{thm:topick}:} We suppose that there is
no $m \geq 1$ \st $f^m$ has infinitely many fixed points. Then
every fixed point for $f^m$ is isolated for any $m$. So we have
the result of Theorem~\ref{thm:holo}. In particular, this means
that $L(f^m) \in \NNN \cup \{0\}$ for all $m \in \NNN$.

In the case that $H_k(M)=0$ for $k>0$, $f_{\ast 0} = (1)$,
multiplication by 1; and for $k>0$, $f_{\ast k}$ is the zero map.
So we have $L(f^m)=1$ for all $m \geq 1$. Therefore, by the
Lefschetz Fixed Point Theorem, $f$ has a fixed point. Furthermore,
by Theorem~\ref{thm:holo}, $\#{\rm Fix}(f^m) \leq 1$ for all $m
\geq 1$.  Therefore, ${\rm Per}(f) =\{1\}$.

When $n \geq 1$,  $f_{\ast 1}$ is an $n \times n$ matrix with
integer entries and with eigenvalues $\lambda_1, \ldots,
\lambda_n$, see \cite{vick}.  So ${\rm trace}(f_{\ast 1}) =
\lambda_1 + \cdots + \lambda_n$ and $L(f^m) = 1- \lambda_1^m -
\cdots - \lambda_n^m$.

We reorder the eigenvalues in order of decreasing modulus (when
two eigenvalues have the same modulus, any choice of order
suffices).

{\bf Case 1:} Suppose that $|\lambda_1|> 1$.

Suppose that $|\lambda_1| = \cdots = |\lambda_k|$ for some maximal
$1 \leq k < n$ (the case  $k=n$ follows similarly). Consider
$\mu_j = \frac{\lambda_j}{|\lambda_1|}$ for $1 \leq j \leq k$.  If
we let $\epsilon>0$ as in Lemma~\ref{lem:modone} we have some
sequence $m_l \to \infty$ \st $Re(\mu_1^{m_l} + \cdots +
\mu_k^{m_l})> k(1-\epsilon)$. Therefore,
$$Re(\lambda_1^{m_l} + \cdots + \lambda_k^{m_l})>
|\lambda_1|^{m_l}k(1-\epsilon) \ {\rm and} \
|Re(\lambda_{k+1}^{m_l} + \cdots + \lambda_n^{m_l})| \leq
(n-k)|\lambda_{k+1}|^{m_l}.$$ Since $|\lambda_1|
>|\lambda_{k+1}|$ there must exist some large enough $m_l$ \st
$$Re(\lambda_1^{m_l} + \cdots + \lambda_k^{m_l})
>1+ |Re(\lambda_{k+1}^{m_l} + \cdots + \lambda_n^{m_l})|.$$ Thus
$1-(\lambda_1^{m_l} + \cdots + \lambda_n^{m_l})$ is negative,
which is not possible by Theorem~\ref{thm:holo}.

{\bf Case 2:} Suppose that $0<|\lambda_1|<1$.

If $\lambda_2, \ldots, \lambda_n=0$ then $L(f)$ cannot be an
integer so this is not possible by the definition of Lefschetz
number. If $|\lambda_2|>0$ we will again show that there exists
some $m \geq 1$ \st $L(f^m)$ is not an integer.  We suppose that
$|\lambda_1| =\cdots = |\lambda_k|$ for some maximal $1 \leq k <
n$ (the case $k=n$ follows similarly). As in case 1, we can choose
$\epsilon>0$ and let $\mu_j = \frac{\lambda_j}{|\lambda_1|}$ for
$1 \leq j \leq k$. From Lemma~\ref{lem:modone}, we have a sequence
$m_l \to \infty$ \st $Re(\mu_1^{m_l} + \cdots + \mu_k^{m_l})>
k(1-\epsilon)$ and so
$$Re(\lambda_1^{m_l} + \cdots + \lambda_k^{m_l})>
|\lambda_1|^{m_l}k(1-\epsilon)\ {\rm and} \
|Re(\lambda_{k+1}^{m_l} + \cdots + \lambda_n^{m_l})| \leq
(n-k)|\lambda_{k+1}|^{m_l}.$$ Since $|\lambda_1|
>|\lambda_{k+1}|$ there must exist some large enough $m_l$ \st
$$|Re(\lambda_1^{m_l} + \cdots + \lambda_k^{m_l})|
>|Re(\lambda_{k+1}^{m_l} + \cdots + \lambda_n^{m_l})|,$$ so
$Re(\lambda_1^{m_l} + \cdots + \lambda_k^{m_l}) \neq 0$.  Since
for large $m_l$, we have $|\lambda_1^{m_l} + \cdots +
\lambda_n^{m_l}|<1$ we know that $1-(\lambda_1^{m_l} + \cdots +
\lambda_n^{m_l})$ cannot be an integer.  This is not possible.

{\bf Case 3:} Suppose that $|\lambda_1| = 1$.  We will show that
the only possibility is that $\lambda =1$, $\lambda_2, \ldots,
\lambda_n=0$ and we have no fixed points.

Suppose that $|\lambda_1|=\cdots =|\lambda_k| =1$ for some maximal
$1 \leq k \leq n$.

{\bf Case 3a:} First suppose that $k \geq 2$.  For large $m$,
$Re(\lambda_{k+ 1}^m + \cdots + \lambda_n^m)$ is very small.
Again, by Lemma~\ref{lem:modone}, for any $\epsilon>0$ we have
some sequence $m_l \to \infty$ \st

$$Re(\lambda_1^{m_l} + \cdots + \lambda_k^{m_l})> k(1-\epsilon).$$

So, since $k \geq 2$, for $m_l$ large enough we have
$1-\lambda_1^{m_l} + \cdots + \lambda_n^{m_l}$ negative, a
contradiction by Theorem~\ref{thm:holo}.

{\bf Case 3b:} Suppose that $k=1$. If $\lambda_1$ were not real
then $\overline{\lambda_1}$ would also be an eigenvalue, so $k$
could not be 1.  Therefore, $\lambda_1 = \pm 1$. We may suppose
that $\lambda_i =0$ for $1 < i \leq n$ (the case where $0
<|\lambda_2|$ is clear because as in Case 2, for large $m$,
$L(f^m)$ will not be an integer). If $\lambda_1=1$ then $L(f^m)=0$
for all $m$ so there are no fixed points of $f^m$.

If $\lambda_1= -1$ then $L(f)=2$ and $L(f^2) =0$.  By the
Lefschetz Fixed Point Theorem, $L(f)=2$ implies that $f$ has a
fixed point. Therefore, $f^2$ has a fixed point and by
Theorem~\ref{thm:holo}, $L(f^2) \geq 1$: so we have a
contradiction.

{\bf Case 4:}  Suppose that all the eigenvalues are zero then
$L(f^m)=1$ for all $m \in \NNN$.  Therefore, $f$ has a fixed point
and ${\rm Per}(f) =\{1\}$. This completes the proof of the
theorem. {\ \hfill\rule[-1mm]{2mm}{3.2mm}}


\begin{thebibliography}{99}

\bibitem{brown}{R.F. Brown, }{\it
The Lefschetz Fixed Point Theorem, } Scott, Foresman and Co.,
Glenview, IL, 1971.

\bibitem{earlham}{C.J. Earle and R.S. Hamilton, }{\it A fixed point
theorem for holomorphic mappings, } Global Analysis, Proc. Symp.
Pure Math., Vol 16, Amer. Math. Soc., Providence, R.I., 1970,
61--65.

\bibitem{holo}{N. Fagella and J. Llibre, }{\it
Periodic points for holomorphic maps via Lefschetz numbers, }
Trans. Amer. Math. Soc. {\bf 352} (2000), 4711--4730.

\bibitem{harris}{L.A. Harris, }{\it Fixed points of holomorphic
mappings for domains in Banach Spaces, }  Abstr. Appl. Anal. {\bf
5 } (2003), 261--274.

\bibitem{hou}{J.X. Hou, }{\it Periodic points of analytical dynamical
systems on hyperbolic Riemann surfaces, } Acta Sci. Natur. Univ.
Sunyatseni {\bf 1} (1988), 28--33.

\bibitem{milnonedim}{J. Milnor, }{\it
Dynamics in one complex variable, } 2nd Ed., Vieweg, 2000.

\bibitem{milngen}{J. Milnor, }{\it Notes on dynamical systems, }
www.math.sunysb.edu/$\sim$jack/DYNOTES/.

\bibitem{nz}{I. Niven and H.S. Zuckerman, }{\it
An introduction to the theory of numbers, }4th Ed., Wiley, NY
1980.

\bibitem{pick}{G. Pick, }{\it \"{U}ber einer eigenschaft der konformen
abbildung kreisf\"{o}rmiger berieche, }  Math. Annalen {\bf 77}
(1916), 1--6.

\bibitem{vick}{J.W. Vick, }{\it Homology Theory, }
2nd Ed., Springer-Verlag, 1994.

\end{thebibliography}
\end{document}